\documentclass[11pt,oneside]{article}
\usepackage{amssymb,amsmath}
\usepackage{url}
\usepackage[english]{babel}
\renewcommand{\le}{\leqslant}
\renewcommand{\ge}{\geqslant}

\newcommand{\RR}{\mathbb{R}}
\newcommand{\ZZ}{\mathbb{Z}}

\newcommand{\NN}{\mathbb{N}}

\newcommand{\FFF}{\mathcal{F}}
\newcommand{\SSS}{\mathcal{S}}
\newcommand{\HHH}{\mathcal{H}}

\newcommand{\MMM}{\mathcal{M}}

\newcommand{\vx}{\mathbf{x}}
\newcommand{\vu}{\mathbf{u}}

\newcommand{\vn}{\mathbf{n}}
\newcommand{\va}{\mathbf{a}}
\newcommand{\vb}{\mathbf{b}}

\newtheorem{lemma}{Lemma}
\newtheorem{theorem}{Theorem}
\newtheorem{proposition}{Proposition}
\newtheorem{problem}{Problem}

\newenvironment{proof}{\textsc{Proof. }}{\ \newline\hspace*{\fill}$\boxtimes$}

\newcommand{\Sl}{\mathrm{SL}}

\begin{document}

\title{Generating random factorisations of polynomial values}

\author{
 Dmitry Badziahin
}

\maketitle

\begin{abstract}
We construct algorithms that efficiently generate random
factorisations of values $P(n)$ as products of two integers, where
$P\in\ZZ[x]$ is a given quadratic or cubic monic polynomial. In
other words, the algorithms produce random triples
$(n,d_1,d_2)\in\ZZ^3$ that solve the Diophantine equation $P(n) =
d_1d_2$. In the case where $P$ is cubic, such an algorithm allows
the construction of an RSA key of $k$ bits that can be described
using about $k/3$ bits of information. We also show how to construct
a solution $(n,d_1,d_2)$ with the ratio $d_1/d_2$ arbitrarily close
to any given positive real number. This proves that among all
solutions $(n,d_1,d_2)$ of $P(n) = d_1d_2$ the ratios $d_1/d_2$ are
dense in $(0,+\infty)$.

\end{abstract}

{\footnotesize{{\em Keywords}: integer factorisation, semiprimes,
RSA key generation

Math Subject Classification 2020: 11Y05, 11T71, 11D99, 94A60}}

\section{Introduction}

An important part of cryptosystems such as RSA is the construction
of a number $N\in \NN$ that is the product of two carefully chosen
prime numbers, $N=pq$. To ensure security, such numbers must satisfy
several conditions, the most important being that $p$ and $q$ are of
similar (and very large) size. For example, in practical RSA
implementations both primes are now typically required to have close
to 1024 bits or more.


Motivated by certain cryptographic applications, the following
natural question arises:
\begin{problem}
Generate an RSA modulus $N$ that has $k$ bits and can be described
using as few bits of information as possible.
\end{problem}
Various authors have addressed this problem by fixing in advance as
many bits of $N$ as possible, see for example~\cite{lenstra_1998,
gra_shp_2008, meng_2013} and references therein. In particular,
in~\cite{lenstra_1998} a heuristic algorithm is presented that
generates $N$ with at most $k/2$ prescribed digits.
In~\cite{gra_shp_2008, meng_2013} different algorithms are proposed
that require slightly fewer bits to be fixed, but whose expected
complexity is proven to be polynomial time (unconditionally
for~\cite{gra_shp_2008} and under GRH for~\cite{meng_2013}). In all
these cases, the asymptotic number of bits required to describe the
output $N$ is $k/2$.

We approach Problem~A from a different perspective, which is also of
independent mathematical interest. Specifically, we construct
factorisations $N=pq$, with $p,q\in\NN$ where in addition $N$ is the
value of a given polynomial $P\in\ZZ[x]$ at an integer argument. In
other words, we aim to produce a randomly chosen integer $n$ such
that $P(n)$ is the product of two large primes of roughly equal
size. If $\deg(P)=d$ then describing the RSA modulus $N = P(n)$
requires specifying only $n$, which has approximately $k/d$ bits. In
this paper, we present algorithms for all quadratic and cubic monic
polynomials, potentially reducing the number of bits needed to
describe $N$ to $k/3$. On the other hand, numbers of the form $P(n)$
might admit more efficient factorisation algorithms, for example,
certain variants of SNFS. Further investigation in this direction is
needed.


The core idea of the method is to rapidly generate random integer
triples $(n,d_1,d_2)$ satisfying the equation
\begin{equation}\label{main_eq}
P(n) = d_1d_2.
\end{equation}
One can then generate many such triples until a solution is found in
which both $d_1$ and $d_2$ are prime. For certain polynomials,
however, such an algorithm will never terminate or will only produce
trivial factorisations. For example, this occurs if $P$ is reducible
over $\ZZ$ or if $P(n)$ is is divisible by a fixed small prime for
all $n\in\ZZ$. For the remaining polynomials, even proving that
$P(n)$ assumes infinitely many semiprime values is a difficult and
in general unsolved problem, see~\cite{ivanec_1978, lemke_2012} for
partial progress. We do not attempt to address this question here,
but note that in practice our algorithms perform well and produce
triples $(n,d_1,d_2)$ with prime $d_1,d_2$ in reasonable time, see
Secion~\ref{sec6} for more details.

For convenience, we restrict our attention to monic polynomials.
With enough care, however, the algorithms can be generalised to
quadratic and cubic polynomials of arbitrary form. Indeed, it is not
hard to observe that for any given $P\in \ZZ[x]$ there exist a monic
polynomial $P^*\in\ZZ[x]$ and constants $a,k\in \ZZ$ such that
$$
aP(x) = P^*(kx).
$$
Therefore all factorisations of $P(n)$ can be obtained from the
factorisations of $P^*(n)$ subject to the additional condition that
$k\mid n$.

In this paper, we present a method for randomly generating triples
$(n,d_1,d_2)$ from the entire set of factorisations of $P(n)$ for
polynomials of the form $P(n) = n^2+c$, $c\in\ZZ$. Although the
method relies only on elementary arguments, to the author was
surprised to find no such result in the literature. For other
quadratic and cubic polynomials, our approach generates triples from
large subsets of all possible factorisations, though not from the
full set. Nevertheless, these subsets are sufficiently large to
establish the following result:

\begin{theorem}\label{th1}
Let $P\in\ZZ[x]$ be a monic polynomial of degree 2 or 3 and
$\alpha\in\RR$, $\alpha\ge 0$. For all $\epsilon>0$ there exist $n,
d_1,d_2\in\ZZ$, which can be computed explicitly, such that
$$
P(n) = d_1d_2\quad\mbox{and}\quad
\left|\frac{d_1}{d_2}-\alpha\right|<\epsilon.
$$
\end{theorem}

{\bf Remark.} It is an easy exercise to verify Theorem~\ref{th1} for
linear polynomials hence we can make this theorem a bit more general
by requiring $P$ to be non-constant and have degree at most 3.

After presenting the necessary background from algebraic number
theory, we begin with a general setup in Section~\ref{sec3}. This
setup applies to polynomials $P$ of arbitrary degree $d$ and
transforms equation~\eqref{main_eq} into another
equation~\eqref{eq3} of degree $d-1$ in $d$ variables. In the
following two sections, we provide algorithms to generate solutions
of this equation for quadratic and cubic polynomials, respectively.
Section~~\ref{sec6} discusses the results of practical
implementations of these algorithms. Finally, Section~\ref{sec7}
presents a slightly different algorithm that enables generation of
all factorisations of the form $n^2+c = d_1d_2$, whereas the
algorithms from the previous sections do not achieve this for some,
and in fact for almost all, polynomials~$P$.

%

\section{Preliminaries and a bit of algebraic number theory}

Let $P\in\ZZ[x]$ have the roots $x_1,x_2,\ldots, x_d$, i.e.
$$
P(x) = (x-x_1)\cdots (x-x_d).
$$
We pick one of these roots and call it $\xi$ and consider the ring
$\ZZ[P]:= \ZZ[x]/ (P(x)\ZZ)$. Notice that for irreducible
polynomials $P$ it is isomorphic to the classical ring $\ZZ[\xi]$.
Also notice that $\ZZ[P]$ is an integer module of dimension $d= \deg
P$. Hence every element $\va\in \ZZ[P]$ can be parameterised by a
$d$-dimensional integer vector $(a_0,\ldots, a_{d-1})\in\ZZ^d$ which
we will also call $\va$.

Given $\va\in\ZZ[P]$, we denote by $P_\va$ the polynomial $P\va(x):=
a_0 + a_1x + \ldots +a_{d-1}x^{d-1}$. for $\va\in\ZZ[P]$ we define
the norm of $\va$ by
$$
N(\va):= \prod_{i=1}^d P_\va(x_i).
$$
One can check that this norm always takes integer values and is
multiplicative, i.e. $N(\va_1\va_2) = N(\va_1)N(\va_2)$. Notice that
for irreducible $P$ this is the classical norm of an algebraic
number $\va\in\ZZ[\xi]$. Finally, observe that for $n\in\ZZ$,
$$
N(n-x) = P(n).
$$
The last observation suggests an idea for generating factorisations
of numbers of the form $P(n)$. Instead, one can try to generate
factorisations of the elements $n-x\in\ZZ[P]$ as a product $\va\vu$,
where $\va,\vu\in\ZZ[P]$.

\section{General setup and a bit of linear algebra}\label{sec3}

Let $P\in\ZZ[x]$ be
\begin{equation}\label{eq1}
P(x) = x^d + c_{d-1}x^{d-1} + \ldots + c_1x+c_0
\end{equation}
and $\va,\vu\in\ZZ[P]$. To compute $\va\vu$, one can multiply two
polynomials $P_\va(x)P_\vb(x)$ and then use~\eqref{eq1} several
times to represent $x^d, x^{d+1}, \ldots, x^{2d-2}$ as linear
combinations of $1,x,\ldots, x^{d-1}$. As the result, we end up with
$$
\va\vu = e_0 + e_1x+\ldots, e_{d-1}x^{d-1},
$$
where $e_i$ for $0\le i<d$ are of the form $e_i = E_{P,i}\va\cdot
\vu$, and $E_{P,i}$ is an $d\times d$ matrix whose entries only
depend on $P$ but not on $\va$ or $\vu$.
By equating $\va\vu=n-x$ and comparing the coefficients for each
power of $x$ on both sides, we derive that $E_{P,i}\va \cdot \vu =
0$ for all $2\le i\le d-1$. In other words, $\vu$ belongs to the
subspace $\langle E_{P, 2}\va, E_{p,3}\va, \ldots,
E_{P,d-1}\va\rangle^\bot$. All vectors in such a subspace can be
parameterised with help of an auxiliary vector $\vx$ in the
following way:
\begin{equation}\label{eq2}
\vu = \vx\wedge E_{P,2}\va \wedge \cdots \wedge E_{P,n-1}\va.
\end{equation}
In fact, since the subspace is two-dimensional, only two coordinates
in $\vx$ suffice to cover it all, but it is more convenient to
consider an arbitrary $\vx$, at least for now. The vector $\vu$ can
be computed in the following way. We first compose an $(d-2)\times
d$ matrix $M$ from the vectors $E_{P,2}\va, \ldots, E_{P,d-1}\va$ as
rows. By $M_{i,j}$, $0\le i\neq j<d$, we denote its minor after
removing the $i'th$ and $j$'th columns from $M$. An important
observation is that $M_{i,j}$ are homogeneous polynomials of degree
$d-2$ in $\va_0, \ldots, \va_{d-1}$. Then we have
\begin{equation}\label{eq4}
u_i = \sum_{j=0\atop j\neq i}^{d-1} (-1)^{\tau(i,j)} x_j M_{i,j},
\end{equation}
where $\tau(i,j)$ is the parity of the permutation of $0,1,\ldots,
d-1$ where $i$ and $j$ are the first two terms and then the
remaining terms follow in ascending order.

Finally, by comparing the coefficient for $x$ in the equation
$\va\vu=n-x$, we derive that $E_{p,1}\va\cdot \vu = -1$.
Substituting~\eqref{eq2} into this equation gives
$$
-(\vx\wedge E_{P,2}\va\wedge\cdots\wedge E_{P,d-1}\va)\cdot
E_{P,1}\va = \vx\wedge E_{P,1}\va\wedge \cdots\wedge E_{P,d-1}\va =
1.
$$
The last equation can be rewritten to
\begin{equation}\label{eq3}
\sum_{i=0}^{d-1}(-1)^ix_i M_i = 1,
\end{equation}
where $M_i$ is constructed similarly to $M_{i,j}$: we first compose
the matrix $M^*$ from the vectors $E_{P,1}\va, \ldots, E_{P,d-1}\va$
as rows and then $M_i$ is its minor after removing $i$'th column.
Notice that $M_i$ are homogeneous polynomials in $a_0, \ldots,
a_{d-1}$ of degree $d-1$.

Now, one can try to generate factorisations $n-x = \vu\va$ in the
following way. Pick a (random) primitive vector $\vx$ and (randomly)
choose a solution $\va$ of the equation~\eqref{eq3}. Then $\vu$ is
computed via~\eqref{eq4}. However, this approach has a couple of
problematic points:
\begin{itemize}
\item We reduce the problem of generating a solution
to~\eqref{main_eq} which is the polynomial equation of degree $d$ to
another polynomial equation~\eqref{eq3} of degree $d-1$. While
reducing the degree is good, it also increases the number of
variables, so finding the general solution of the new equation may
become tricky. On the other hand, for $d\le 3$, we end up with
polynomials of degree 2 which are well understood.

\item From~\eqref{eq4} we get that $\vu$ is a polynomial of $\va$
of degree $d-2$. Going back to the equation $P(n) = N(\va)N(\vu)$,
we get that $N(\va)$ is a polynomial of degree $d$ of $\va$, while
$N(\vu)$ is one of degree $d(d-2)$. Hence for $d\ge 4$ and generic
solutions $\va, \vu$ the values of $N(\va)$ and $N(\vu)$ will differ
be too much and hence will not be cryptographically secure.
\end{itemize}

Also notice that the described approach does not always cover all
possible triples $(n,d_1,d_2)$ that satisfy~\eqref{main_eq}. For
example, for $P(x) = x^2+5$ we have $1^2 + 5 = 6 = 2\cdot 3$, but
there are no $\va,\vu\in\ZZ[P]$ such that $N(\va) = 2$ and
$N(\vu)=3$. In the last section we will provide the method that
generates all possible triples for quadratic polynomials of the form
$P(x) = x^2+c$. With a bit of an effort, one can generalise it to
all quadratic polynomials. We leave that task to an enthusiastic
reader.

\section{Quadratic Polynomials}\label{sec4}

Let $P$ be a quadratic polynomial:
$$
P(x) = x^2 + c_1x+c_0.
$$
Then we have $\va = a_0 +a_1x$ and $\vu = u_0 + u_1x$ and all the machinery from the previous section drastically simplifies. In particular, the equation
$$
\va\vu = n-x
$$
gives only one equation which comes from comparing the coefficients
for $x$ on both sides of the equation $$\va\vu = a_0u_0 - c_0a_1u_1
+ (a_1u_0 + (a_0 - c_1a_1)u_1)x.$$ This gives the equation
$$
a_1u_0 + (a_0 - c_1a_1)u_1=-1.
$$
It can be rewritten in the form
\begin{equation}\label{eq5}
\det\left(\begin{array}{cc}
c_1a_1-a_0& a_1\\
u_0&u_1\\
\end{array}\right) = 1.
\end{equation}
Notice that the solutions of this equation are in one-to-one
correspondence with the elements in $\Sl_2(\ZZ)$. But this set is
well understood, so we can provide an algorithm that generates its
elements.

Let $A = \left(\begin{array}{cc}
a_0^*& a_1\\
u_0&u_1\\
\end{array}\right)$. Without loss of generality we may assume
that $|u_1|\le |a_1|$ because otherwise we just swap $\va$ and
$\vu$. Next, we observe that $\gcd(a_0^*, a_1)=1$ otherwise the
determinant of the matrix can not be 1. For a fixed $(a_0^*, a_1)$
it is not hard to check that there are only two possibilities for
the pair $(u_0, u_1)$ to satisfy~\eqref{eq5}. They can be obtained
in the following way. First, one computes the continued fraction of
$$
\frac{\pm a_0^*}{a_1} = [r_0;r_1,r_2,\ldots, r_s].
$$
Then for $(s-1)$'st convergent $p_{s-1}/q_{s-1}$ of this fraction we
have $(\pm a_0^* q_{s-1}) - a_1p_{s-1} = (-1)^{s-1}$. Finally, an
element of $\Sl_2(\ZZ)$ is derived by taking $u_0 = (-1)^{s-1}
p_{s-1}, u_1 = \pm(-1)^{s-1}q_{s-1}$. By choosing the sign in $\pm
a_0^*$, we end up with two elements of $\Sl_2(\ZZ)$.

After generating the matrix $A$, one computes $a_0 = c_1a_1-a_0^*$
and from the equation $\va\vu=n-x$,
$$
n = (c_1a_1-a_0^*)u_0 - c_0a_1u_1.
$$
Taking norms on both sides gives
$$
(a_0^2 - c_1a_0a_1 + c_0a_1^2)(u_0^2 - c_1u_0u_1 + c_0u_1^2) = P(n).
$$
That is our random factorisation of the form $P(n) = d_1d_2$.

Rough estimates on $d_1,d_2$ give $d_1d_2\lesssim ||\va||^2_\infty
||\vu||^2_\infty ||(c_0, c_1)||_\infty^2$. If one needs $P(n)$ to be
of certain size, e.g. $|P(n)|\asymp M$ where $M$ is for example
$2^{2048}$, then one can require $\va$ to be of the size
$$
||\va||_\infty \le \left(\frac{M}{c^2}\right)^{1/4}, \quad\mbox{where}\quad c:= ||(c_0, c_1)||_\infty.
$$

Now we are ready to sketch the algorithm for generating random
factorisations~\eqref{main_eq}.
\begin{enumerate}
\item Randomly generate integers $a^*_0, a_1$ such that $|c_1a_1-a_0^*|, |a_1| \le \sqrt[4]{M/c^2}$ and $\gcd(a_0^*, a_1)=1$.
\item Randomly choose the sign of $\pm \frac{a_0^*}{a_1}$, compute
its continued fraction expansion and finally compute $p_{s-1}$ and
$q_{s-1}$ for its $(s-1)$'st convergent.
\item Compute $u_0 = (-1)^{s-1}p_{s-1}$, $u_1 = \pm (-1)^{s-1}q_{s-1}$.
\item Compute $a_0 = c_1a_1-a_0^*$ and $n = a_0u_0 - c_0a_1u_1$.
\item Compute $d_1=a_0^2 - c_1a_0a_1 + c_0a_1^2$, $d_2 = u_0^2 - c_1u_0u_1+c_0u_1^2$.
\end{enumerate}

Usually, in cryptographic applications we need $d_1$ and $d_2$ to be prime. In that case, we add one more step:
\begin{enumerate}
\item[6.] Check if $d_1$ and $d_2$ are prime. If at least one of them is not then go back to Step~1.
\end{enumerate}

The standard properties of continued fractions imply that if
$a_0^*/a_1 = [r_0;r_1,\ldots, r_s]$ then $a_0^*/u_0 = [r_s; r_{s-1},
\ldots, r_0]$. Therefore, for all pairs $(\beta,\gamma)\in\RR^2$ one
can fix several first and last partial quotients in the fraction
$a_0^*/a_1$ to construct matrix $A$ such that $a_0^*/a_1$ and
$u_0/u_1$ are as close to $\gamma$ as we wish and $a_0^*/u_0$ and
$a_1/u_1$ are as close to $\beta$ as we wish. In particular, we can
construct a sequence of matrices
$\begin{bmatrix}a^*_{0,k}&a_{1,k}\\u_{0,k}&u_{1,k}\end{bmatrix}$
such that
\begin{equation}\label{eq11}
\lim_{k\to\infty} \frac{a_{0,k}^*}{a_{1,k}} = \lim_{k\to\infty}
\frac{u_{0,k}}{u_{1,k}} = \gamma,\quad \lim_{k\to\infty}
\frac{a_{0,k}^*}{u_{0,k}} = \lim_{k\to\infty}\frac{a_{1,k}}{u_{1,k}}
= \beta.
\end{equation}
By applying the formulae from Steps~4 and~5 of the algorithm, we get
$$
a_{0,k}/u_{1,k} \to (c_1-\gamma)\beta,\quad n_k/u_{1,k}^2 \to
-(\gamma^2-c_1\gamma+c_0)\beta,
$$$$
d_{1,k}/u_{1,k}^2\to (\gamma^2-c_1\gamma+c_0)\beta^2,\quad
d_{2,k}/u_{1,k}^2 \to \gamma^2-c_1\gamma+c_0
$$
as $k$ tends to infinity. Finally, we derive
$$
d_{1,k}/d_{2,k} \to \beta^2, \quad n_k/d_{2,k}\to -\beta.
$$
As $\beta\in\RR$ is arbitrary, this verifies Theorem~\ref{th1} for
quadratic polynomials $P$. It also provides the method to generate
multiplications~\eqref{main_eq} with $\big|\frac{d_1}{d_2} -
\alpha\big|<\epsilon$ for arbitrary $\alpha\ge 0$ and arbitrarily
small $\epsilon$. One just needs to construct the matrix $A$ with
$a_0^*/u_0$ close enough to $\sqrt{\alpha}$.

For the purposes of the next section, we need to make the last
statement stronger. Fix $q\in\NN$ and let $c_0 = qc_0^*$. For any
$\epsilon>0$ we can always find a continued fraction $[r_k;r_{k-1},
\ldots, r_0]$ such that $\big|[r_k;r_{k-1},\ldots, r_0] -
\beta\big|\le \epsilon$ and its denominator is an integer multiple
of $q$. We can also find a continued fraction $[l_0;l_1,\ldots,
l_k]$ such that $\big|[ l_0;l_1,\ldots, l_k] - \gamma|\le \epsilon$.
Next, there exists $l\in\NN$ and a continued fraction
$[m_l;m_{l-1},\ldots, m_0]$ such that
$$
A_{l_0}A_{l_1}\cdots A_{l_k}A_{m_0}\cdots A_{m_l}\equiv
\mathrm{Id}\;(\mathrm{mod }\; q).
$$
That all shows that we can construct a sequence of matrices
$\begin{bmatrix}a^*_{0,k}&a_{1,k}\\u_{0,k}&u_{1,k}\end{bmatrix}$
such that on top of~\eqref{eq11} they satisfy $q\mid u_{0,k}$ or
equivalently $u_{0,k} = qu^*_{0,k}$.

Notice that
$$
n_k = a_{0,k}u_{0,k} - c_0a_{1,k}u_{1,k} = q (a_{0,k}u^*_{0,k} -
c_0^* a_{1,k}u_{1,k}) = qn^*_k
$$
and then
$$
P(n_k) = qP^*(n^*_k), \quad\mbox{where}\quad P^*(x) = qx^2 + c_1x +
c_0^*;
$$
$$
d_{1,k} = a_{0,k}^2 - c_1a_{0,k}a_{1,k} + qc_0^*a_{1,k}^2,\quad
d_{2,k} = q(qu^{*2}_{0,k} - c_1u_{0,k}^*u_{1,k} + c_0^*u_{1,k}^2) =
qd_{2,k}^*.
$$
In other words, we generate the sequence $(n_k^*, d_{1,k},
d^*_{2,k})$ of solutions of the equation $P^*(n) = d_1d_2$ with
$$
n_k^*/d_{2,k}^* \to -\beta,
$$
where $P^*\in\ZZ[x]$ is not necessarily monic.

\section{Cubic polynomials}\label{sec5}

We continue with a more complicated case of cubic polynomial
$$
P(x) = x^3+c_2x^2+c_1x+c_0.
$$
Then $\va = a_0+a_1x+a_2x^2$, $\vu = u_0+u_1x+u_2x^2$ and
$$
\va\vu = (a_0u_0-c_0a_2u_1+c_0(c_2a_2-a_1)u_2) +
(a_1u_0+(a_0-c_1a_2)u_1+((c_1c_2-c_0)a_2-c_1a_1)u_2)x
$$$$
+ (a_2u_0 + (a_1-c_2a_2)u_1+ (a_0-c_2a_1+(c_2^2-c_1)a_2)u_2)x^2.
$$
Then we compute
$$
E_{P,1} = \begin{pmatrix} 0&1&0\\
1&0&-c_1\\
0&-c_1&c_1c_2-c_0
\end{pmatrix},\qquad E_{P,2} = \begin{pmatrix}
0&0&1\\
0&1&-c_2\\
1&-c_2&c_2^2-c_1
\end{pmatrix}
$$
and
$$M^* = \begin{pmatrix}
a_1&a_0-c_1a_2&-c_1a_1+(c_1c_2-c_0)a_2\\
a_2&a_1-c_2a_2&a_0-c_2a_1+(c_2^2-c_1)a_2
\end{pmatrix}.
$$
Next, we compute
$$
M_0 = a_0^2-c_2a_0a_1+c_1a_1^2+(c_2^2-2c_1)a_0a_2-(c_1c_2-c_0)a_1a_2
+ (c_1^2-c_0c_2)a_2^2,
$$
$$
M_1 = a_0a_1-c_2a_1^2+c_2^2a_1a_2 - (c_1c_2-c_0)a_2^2,\quad M_2 =
a_1^2-a_0a_2-c_2a_1a_2+c_1a_2^2.
$$

For simplicity, in the equation~\eqref{eq3} we set $x_0=0$. Under
this condition we can not guarantee anymore that each solution
$\va\vu = n-x$ can be obtained. But since the dimension of $\langle
E_{P,2}\va\rangle^\bot$ is two-dimensional, not too many of them
will be lost. For simplicity, denote $(x_1,x_2) = (-x,y)$.
Then~\eqref{eq3} can be rewritten as
\begin{equation}\label{eq6}
a_0(xa_1-ya_2) + Aa_1^2 + Ba_1a_2 + Ca_2^2 = 1,
\end{equation}
where $A = y-c_2x, B = -c_2(y-c_2x)$ and $C = c_1y - (c_1c_2-c_0)x$.

A solution of~\eqref{eq6} can be generated in the following way.
First, notice that $x$ and $y$ must be coprime as otherwise the
equation~\eqref{eq3} does not have solutions. Hence, as the first
step, we randomly choose two coprime numbers $(x,y)$. Then we notice
$$
Aa_1^2 + Ba_1a_2 + Ca_2^2 \equiv 1\;(\mathrm{mod}\; xa_1-ya_2).
$$
denote the modulus $xa_1-ya_2$ by $Q$. Assuming that $y$ and $Q$ are
coprime, we derive $a_2 \equiv \frac{x}{y}a_1$ (mod $Q$) and then
the congruence transforms to
\begin{equation}\label{eq7}
\left(A + \frac{Bx}{y} + \frac{Cx^2}{y^2}\right)a_1^2 \equiv
1\;(\mathrm{mod}\; Q).
\end{equation}
The last inequality can be quickly solved as soon as the prime
factorisation of $Q$ is known.

In principle, the case $\gcd(Q,y)=d>1$ can also be dealt easily. We
then must have $a_1\equiv 0$ (mod $d$) and the congruence transforms
to $Ca_2^2\equiv 1$ (mod $d$) which can be solved again. For the
remaining part $Q/d$ of the modulus we have $\gcd(y/d, Q/d)=1$ and
we can proceed in the usual way.

As the next step of the algorithm, we choose random $Q$. It either
needs to be not large enough so that one can quickly factorise it.
Or one can require $Q$ to be prime. In both cases, some of the
factorisations will not be covered. We will stick with the second
option and select a random prime $Q$. If~\eqref{eq7} has no
solutions modulo that number, we pick another $Q$. Since half of all
numbers between $1$ and $Q-1$ are quadratic residues, we expect the
congruence~\eqref{eq7} has solutions for approximately half of
primes $Q$, therefore we should be able to find an appropriate $Q$
quickly.

Now suppose that for a given prime $Q$ the solution of~\eqref{eq7}
is $a_1\equiv \pm a$ (mod $Q$). We can also compute $a_2\equiv
\frac{x}{y}a_1 \equiv \pm b$ (mod $Q$). Then by construction we have
$xa_1-ya_2$ and $Aa_1^2 +Ba_1a_2 + Ca_2^2 - 1$ are multiples of $Q$.

Next, from the general solution $a_1 = \pm a + kQ$, $a_2 = \pm b +
lQ$ we select $k,l\in\ZZ$ such that $xa_1 - ya_2 = Q$. That can be
done quickly, since
$$
xa_1 - ya_2 = \pm (xa- yb) + Q(xk-yl).
$$
and we can quickly find $k,l\in \ZZ$ such that $xk-yl = 1 -
\frac{xa-yb}{Q}$.

Next, we find $a_0$ from equation~\eqref{eq6}:
\begin{equation}\label{eq9}
a_0 = \frac{1 - Aa_1^2-Ba_1a_2 - Ca_2^2}{Q}.
\end{equation}
By construction of $a_1,a_2$, this number is always integer.

Then we compute $\vu$ from the formula~\eqref{eq4}:
\begin{equation}\label{eq8}
u_0 =-(a_0 - c_2a_1+(c_2^2-c_1)a_2)x - (a_1-c_2a_2)y,\quad u_1 =
a_2y,\quad u_2 = a_2x
\end{equation}
and $n$ from the product $\va\vu$:
\begin{equation}\label{eq10}
n = a_0u_0-c_0a_2u_1+c_0(c_2a_2-a_1)u_2.
\end{equation}
Finally, taking norms gives us the factorisation $P(n) =
N(\va)N(\vu)$.

To conclude, we sketch the algorithm for generating the solutions
of~\eqref{main_eq} with cubic $P(x)$.

\begin{enumerate}
\item Generate random numbers $x,y$ such that $|x|, |y|\le B_1$ and
$\gcd(x,y)=1$. Here $B_1$ is a given upper bound.

\item Generate random prime $Q$ such that $Q\le B_2$, $\gcd(Q,y)=1$ and the
congruence~\eqref{eq7} has a solution $a_1$. Here $B_2$ is another
given upper bound.

\item Solve~\eqref{eq7} and let $a$ be one of its solutions. Compute $b\equiv
\frac{x}{y}a$ (mod $Q$).

\item Compute $a_1=a+kQ$ and $a_2=b+lQ$ where $k,l\in\ZZ$ are such that
$xk-yl=1-\frac{xa-yb}{Q}$. These values can be computed with help of
the extended Euclidean algorithm.

\item Compute $a_0$ by~\eqref{eq9} and $\vu$ by~\eqref{eq8}.

\item Finally, compute $n$ by~\eqref{eq10} and $d_1=N(\va)$,
$d_2=N(\vu)$.
\end{enumerate}

As for the case of quadratic polynomials, if we need $d_1,d_2$ to be
prime, we add one more step:

\begin{enumerate}
\item[7.] Check if $d_1$ and $d_2$ are prime. If at least one of them is not then go back to Step~1.
\end{enumerate}

Now we prove Theorem~\ref{th1} for cubic polynomials. Notice that
without loss of generality one may assume that $c_0\neq 0$. Indeed,
if $c_0=0$ then we can replace the polynomial $P(x)$ by an
equivalent one $Q(x) = P(x+1)$. Choose the pair of parameters $(x,y)
= (1,c_2)$. Then the equation~\eqref{eq6} simplifies to
$$
c_0a_2^2 -(c_2a_2-a_1)a_0 = 1.
$$
By changing the variable $b:=c_2a_2-a_1$, we derive the
equation~\eqref{main_eq} for a given quadratic polynomial $P(n) =
c_0 n^2 -1$. From the previous section we know that for all
$\beta\in\RR$ one can construct the sequence $a_{0,k}, a_{2,k}, b_k$
of solutions of this equation such that
$$
\lim_{k\to\infty} \frac{b_k}{a_{2,k}} = \lim_{k\to\infty}
\frac{c_0a_{2,k}}{a_{0,k}}= \beta.
$$
This implies $a_{2,k} = \frac{\beta}{c_0} a_{0,k} + o(a_{0,k})$ and
$a_{1,k} = \frac{(c_2-\beta)\beta}{c_0}a_{0,k} + o(a_{0,k})$. For
convenience, we will write these equations as $a_{2,k} \sim
\frac{\beta}{c_0} a_{0,k}$ and $a_{1,k}\sim
\frac{(c_2-\beta)\beta}{c_0}a_{0,k}$.

We also notice that~\eqref{eq8} simplifies to
$$
\vu = (c_1a_2 - a_0, c_2a_2, a_2)
$$
and after substituting $a_{0,k}, a_{1,k}, a_{2,k}$ into it, we
derive
$$
\vu_k \sim \left(\frac{c_1\beta-c_0}{c_0}, \frac{c_2\beta}{c_0},
\frac{\beta}{c_0}\right) a_{0,k}.
$$

Notice that the norm of a cubic irrational $a+bx + cx^2$ equals
$$
N(a+bx+cx^2) = (a+bx+cx^2)(a+bx_2+cx_2^2)(a+bx_3+cx_3^2)
$$$$
=a^3 -c_2a^2b + (c_2^2-2c_1)a^2c + c_1ab^2 + (3c_0 - c_1c_2)abc + (c_1^2 - 2c_0c_2)ac^2 - c_0b^3 + c_0c_2b^2c - c_0c_1bc^2 + c_0^2c^3.
$$
Careful computations then give
$$
N(\va_k) \sim
\frac{(\beta^3-c_2\beta^2+c_1\beta-c_0)^2}{c_0^2}a^3_{0,k}, \quad
N(\vu_k)\sim \frac{\beta^3-c_2\beta^2+c_1\beta-c_0}{c_0} a^3_{0,k}.
$$
Finally,
$$
\lim_{k\to\infty} \frac{N(\va_k)}{N(\vu_k)} =
\frac{\beta^3-c_2\beta^2+c_1\beta-c_0}{c_0}.
$$
As $\beta$ runs through all real numbers, the last expression also
takes every value in $\RR$. That confirms Theorem~\ref{th1} for
cubic polynomials $P$.

\section{Numerical results}\label{sec6}

The algorithms described in the previous sections were implemented
in C++ using the NTL library for large integer arithmetic. They were
tested on a single core of an Intel Core I5-14500 processor, and all
reported times correspond to this setup. Since the algorithms are
probabilistic, their runtimes before producing output can vary
significantly between runs. To mitigate this variability, each
algorithm was executed five times, and the reported times are the
averages over these runs. the implementation prioritised code
clarity over execution speed. With further optimisations and by
utilising all available processor cores, the computations can be
accelerated by a factor of 10 to 30.

{\bf Quadratic polynomials.} Here we choose $P(x) = x^2+1$ and
search for the triples $(n,d_1,d_2)$ such that $P(n) = d_1d_2$ and
both $d_1$ and $d_2$ are prime. We do not provide all the particular
factorisations that were found with the algorithm as they are
random. As a demonstration, we provide one: for $d_1d_2$ bounded by
200 bits, one of the factorisations was
$$
\scriptstyle 10356460241698236473356068520^2 + 1 =
869011498091105945959277606189 \cdot 123423302192753838490702309.
$$
The results of the algorithm for different sizes of $d_1d_2$ are
presented in the following table.

\begin{center}
\noindent\begin{tabular}{|c|c|c|} \multicolumn{3}{c}{{\bf Table 1.}
Results for $P(x) = x^2+1$.}\\
\hline Bound on $d_1d_2$& Time &Number of tries before finding prime $d_1$ and $d_2$\\
\hline 100 bits&9 msec& 318\\
\hline 200 bits&38 msec& 2261\\
\hline 400 bits&147 msec& 4135\\
\hline 800 bits&1.57 sec& 17255\\
\hline 1024 bits&9.86 sec& 68055\\
\hline 2048 bits&66 sec& 118687\\
\hline
\end{tabular}
\end{center}

We have also implemented an algorithm that generates triples
$(n,d_1,d_2)$ with $P(n) = d_1d_2$ such that the continued fraction
of $a_0/u_0$ starts with that of a fixed fraction $p/q$. As we have
shown in Section~\ref{sec4}, this guarantees that the fraction $a/b$
is close to the number $(p/q)^2$. We run that algorithm for
$\frac{p}{q}=\frac{39}{22}$ and $\frac{p}{q}=\frac{8545}{4821}$ that
are convergents to $\sqrt{\pi}$.

For $\frac{p}{q}=\frac{39}{22}$ we set $d_1d_2$ to be up to 10
digits long. Then the algorithm produced
$$
78457^2+ 1 = 6155500850 = 139186 \cdot 44225,\quad
\frac{139186}{44225} \approx 3.147.
$$

For $\frac{p}{q}=\frac{8545}{4821}$ we set $d_1d_2$ to be 20 digits
long and got
$$
\scriptstyle 2184324567^2 + 1\; =\; 4771273813999737490\; =\;
3871614562\; \cdot\; 1232373145,\quad
\frac{3871614562}{1232373145}\; \asymp\; 3.14159277
$$
Since we did not require $d_1$ and $d_2$ to be prime, the algorithm
completed almost instantly even for huge bounds on $d_1d_2$.
Therefore we do not provide the timing here.


\smallskip
{\bf Cubic polynomials. } Here we choose $P(x) = x^3+x+1$ and
generate random triples $(n,d_1,d_2)$ such that $P(n) = d_1d_2$ with
both $d_1$ and $d_2$ being prime. As described in
Section~\ref{sec5}, the algorithm uses the bounds $B_1$ and $B_2$ on
primitive pairs $(x,y)\in\ZZ^2$ and (prime) numbers $Q$
respectively. In the implementation of the algorithm one can tune
them in order to get the most efficient runtime and the most
appropriate output. Also, even for fixed $B_1$ and $B_2$ the outputs
can vary in a wide range. Values $P(n)$ can easily differ by 5 -- 7
digits from output to output. For example, when $Q$ was taken in the
range between 1 and 100 and $||\vx||_\infty$ was bounded by 1000,
the five consecutive outputs of the algorithm were:
$$
\begin{array}{rrl}
\mathbf{1.}& \scriptscriptstyle P(-1651600654673606743145) = & \scriptscriptstyle -134232642580181579009479643127 * 33562708346658329124432524590457447;\\
\mathbf{2.}& \scriptscriptstyle P(-373568507089772975507364265491) = & \scriptscriptstyle -100342782433902695905377286791922956278293\\
&&\scriptscriptstyle  * 519546748147127834470514792753339521946900823977;\\
\mathbf{3.}& \scriptscriptstyle P(14333236722773622209533415524) = &\scriptscriptstyle 130889291060184889633347224892063375187 \\
&&\scriptscriptstyle * 22497212231386855181739997244013731374992645927;\\
\mathbf{4.}& \scriptscriptstyle P(2422442613890813481822026798646190) = &\scriptscriptstyle 6437736445907770427561169007575986475952878237\\
&&\scriptscriptstyle  * 2208143533333497826538693266693240061269645923787593043;\\
\mathbf{5.}& \scriptscriptstyle P(98616711516091188636140618879) = &\scriptscriptstyle 3827966952584582054331533888109835851017\\
&&\scriptscriptstyle  *
250543632316242226506381013500672643566534801607.
\end{array}
$$

Next, notice that $d_1 = N(\va)$, written as a polynomial of $\va$,
has degree $3$, while $d_2 = N(\vu)$ is of degree 3 over $\va$ and
also degree 3 over $\vx$. Therefore, in the generic triple
$(n,d_1,d_2)$ the value $|d_2|$ is approximably $||\vx||^3_\infty$
times bigger than $|d_1|$. If we choose the bound $B_1$ on
$||\vx||_\infty$ too large then the sizes of $d_1$ and $d_2$ in the
triple will become unbalanced. One of them will be much smaller than
the other which may make the number $P(n)$ vulnerable to
factorisation algorithms such as ECM. On the other hand, if we only
consider small vectors $\vx$ we restrict ourselves to a small
proportion of all solutions $(n,d_1,d_2)$ of the
equation~\eqref{main_eq}. In our implementation we chose
$||\vx||_\infty\le 1000$. The results of the algorithm are provided
in the table below. Again, for each set of the parameters we ran it
five times and provide mean results for all the tries.
\begin{center}
\noindent\begin{tabular}{|c|c|c|c|} \multicolumn{4}{c}{{\bf Table
2.} Results for $P(x) = x^3+x+1$.}\\
\hline Bound on $Q$& Time &Number of tries&Size of $d_1d_2$\\
\hline 100&37 msec& 1509&86 digits\\
\hline 50 bits&688 msec& 4908&159 digits\\
\hline 100 bits&11.1 sec& 22194&260 digits\\
\hline 200 bits&141 sec& 59774&440 digits\\
\hline 400 bits&76.45 min& 281770&794 digits\\
\hline
\end{tabular}
\end{center}

\section{Factorizations of numbers $n^2+c$}\label{sec7}

Here we outline the algorithm which ensures that we are generating
the solutions of~\eqref{main_eq} among all possible factorisations
of the values $P(n)$. We start by classifying all the factorizations
of numbers of the form $n^2+1$, $n\in\ZZ$. It is based on the
following fact:

\begin{proposition}\label{th1}
There exists a bijection $\phi$ between the free semigroup $\SSS$
generated by $2\times 2$ matrices:
$$
\SSS:= \langle A_n, n\in \NN\rangle, \quad A_n: = \left(\begin{array}{cc}0&1\\
1&n
\end{array}\right)
$$
and the set $\FFF$ of factorizations of $n^2+1$
$$
\FFF:=\{(n,d_1,d_2)\in \ZZ_{\ge 0}\times\NN^2\;:\; n^2+1 = d_1d_2,\;
d_1<d_2\}
$$
defined as follows:
$$
\phi(M):= (n,d_1,d_2),\quad \text{where}\quad \left(\begin{array}{cc}d_1&n\\
n&d_2
\end{array}\right) = M^tM.
$$
\end{proposition}

\proof

Since for all $n\in\NN$ $\det(A_n)=-1$, we have that $\det(M) = \pm
1$ for every matrix $M\in\SSS$. Notice that $M^tM$ is symmetric and
therefore
$$
n^2 - d_1d_2 = \det(M^tM) = 1\quad\mbox{for some}\quad
n,d_1,d_2\in\NN.
$$
So $n^2+1 = d_1d_2$ and for every $M\in\SSS$, $\phi(M)$ indeed maps
to some factorization of $n^2+1$ .

Now we show that $\phi$ is surjective. For a given factorization
$(n,d_1,d_2)$ we show that its preimage is nonempty. We do it by
induction on $d_1$. Firstly note that for the trivial factorization
$n^2+1 = 1\cdot (n^2+1)$ the preimage contains $A_n$:
$$
\phi(A_n) = (n,1,n^2+1).
$$
Now assume that for all factorizations $(n,d_1,d_2)$ with $d_1\le D$
the preimages of $(n,d_1,d_2)$ are nonemtpty. Let's consider the
factorization $(n,d_1,d_2)\in \FFF$ with the smallest possible
$d_1>D$. One can easily check that in this case $d_1\le n<d_2$.
After dividing $n$ by $d_1$ with the remainder:
$$
n = qd_1+r, \quad r<d_1\le n
$$
we get that $r^2+1$ divides $d_1$ too and therefore in view of
$r<d_1$ the initial factorization transforms to
$$
r^2+1 = d_1^*d_2^*,\quad d_1^*\le r<d_2^* = d_1.
$$
In other words $d_1^*\le D$, so the preimage of $(r,d_1^*,d_2^*)$ is
nonempty, in particular it includes some matrix $M\in\SSS$. Observe
that
$$
\left(\begin{array}{cc}d_1&n\\
n&d_2
\end{array}\right) =  A_q^t \left(\begin{array}{cc}d_1^*&r\\
r&d_2^*
\end{array}\right) A_q
$$
which finishes the proof of surjectivity.

To verify the injectivity of $\phi$ notice that from
$$
\begin{pmatrix}
d_1&n\\
n&d_2
\end{pmatrix} = A_{n_d}^tA_{n_{d-1}}^t\cdots A_{n_1}^tA_{n_1}\cdots
A_{n_d}
$$
the matrix $A_{n_d}$ is uniquely determined by $n_d = \lfloor
\frac{n}{d_1}\rfloor$. Then inductively, the whole matrix $M$ is
uniquely determined from the triple $(n,d_1,d_2)$.
\endproof



Now consider the general case of factorizations
$$
n^2 + c = d_1d_2, \quad d_1\le d_2.
$$
We first focus on the slightly easier case of $c>0$ and prove the
following auxiliary result:

\begin{lemma}\label{lem1}
Let $c\in\NN$ be fixed. There are only finitely many triples
$(n,d_1,d_2)\in \ZZ_{\ge 0}\times\NN^2$ such that $n^2+c = d_1d_2$
and $n< d_1\le d_2$.
\end{lemma}

\proof  Assume that $d_1 = n+k$ with $k\ge 0$. Then
$$
n^2+c = d_1d_2\ge(n+k)^2 >n^2+k^2.
$$
Therefore $k<\sqrt{c}$ and there are only finitely many
possibilities for the parameter $k$. Now we show that for any fixed
parameter $k$ there are only finitely many solutions $(n,d_1,d_2)$.
It would mean that $n+k\mid n^2+c$. Indeed,
$$
n+k \;|\; n^2+c \quad \Rightarrow \quad n+k\;|\; k^2 + c.
$$
Thus $n+k$ divides the fixed number $k^2+c$, but there are only
finitely many possibilities for~$n$ with this condition.
\endproof

Denote by $D_c$ the finite set of all triples
$(n,d_1,d_2)\in\ZZ_{\ge 0}\times \NN^2$ with $n< d_1\le d_2$ and
$n^2+c = d_1d_2$. Also denote by $\HHH_c$ the associated set of
matrices
$$
\HHH_c:= \left\{\left(\begin{array}{cc}d_1&n\\
n&d_2
\end{array}\right)\;:\; (n,d_1, d_2)\in D_c\right\}.
$$
Notice that the proof of Lemma~\ref{lem1} also suggests the
algorithm for computing the set $D_c$.


\begin{proposition}\label{th2}
Let $c$ be a natural number. There exists a bijection between
$\SSS\times \HHH_c$ and the set $\FFF_c$ of triples $(n,d_1,d_2)$:
$$
\FFF_c:=\{(n,d_1, d_2)\in \ZZ_{\ge 0}\times \NN^2\;:\; n^2+c =
d_1d_2, d_1\le d_2\}.
$$
It is given by:
$$
\phi(M,H):= (n,d_1,d_2),\quad \text{where}\quad \left(\begin{array}{cc}d_1&n\\
n&d_2
\end{array}\right) = M^tHM.
$$
\end{proposition}

\proof We have that $\det(M) = \pm 1$ for every matrix $M\in\SSS$.
Also by construction, $\det(H) = c$ and $M^tHM$ is a symmetric
matrix of the form $\begin{bmatrix}d_1&n\\n&d_2\end{bmatrix}$. So
$n^2+c = d_1d_2$ and for every $M\in\SSS, h\in \HHH_c$, $\phi(M)$
indeed maps to some factorization of $n^2+c$.

We show that the preimage of $(n,d_1, d_2)$ under $\phi$ is
nonempty. We do it by induction on $n$. If $d_1>n$ then by
definition $\begin{bmatrix}d_1&n\\n&d_2\end{bmatrix} \in \HHH_c$ and
we take $M = \mathrm{Id}$. Hence we can assume that $d_1\le n$ and
then proceed as in the proof of Proposition~\ref{th1}. We can write
$$
\left(\begin{array}{cc}d_1&n\\
n&d_2
\end{array}\right) =  A_q^t \left(\begin{array}{cc}a^*&r\\
r&d_1
\end{array}\right) A_q\quad\mbox{where}\quad n=qd_1+r,\quad
r<d_1\le n.
$$
Therefore $\phi$ is surjective.

The injectivity of $\phi$ is verified analogously as in the proof of
Proposition~\ref{th1}.

\endproof

The case $c\le 0$ is slightly more difficult. By following the
arguments of Proposition~\ref{th1}, in this case we may end up with
$d_1^*\le 0$. Luckily, the factorisations $n^2+c = d_1d_2$ with
non-positive $d_1$ may only happen for small values of $n$,
therefore the set of such factorisations can be easily computed.

Denote by $N_c$ the set of all triples $(n,d_1,d_2)\in \ZZ_\ge
0\times\ZZ_{\le 0}\times\NN$ that satisfy $n^2+c = d_1d_2$. One can
easily check the following fact:
\begin{lemma}\label{lem2}
If $-c$ is a perfect square then the set $N_c$ is infinite.
Otherwise it is finite.
\end{lemma}
\proof Indeed, if $c=-k^2$ then $k^2-k^2 = 0\cdot d$ for all
$d\in\NN$ gives us an infinite subset of elements $(k,0,d)$ in
$N_c$. Otherwise $n^2+c$ is never equal to zero and therefore for
$(n,d_1,d_2)\in N_c$ one must have $n^2+c<0$. It is satisfied for
finitely many values $n$ and each of those values $n$ gives finitely
many factorisations $n^2+c=d_1d_2$.\endproof

Now given an element $\vn:= (n,d_1,d_2)\in N_c$ we define the set
$\MMM_\vn$ of matrices:
$$
\MMM_\vn:=\left\{\left(\begin{array}{cc}u&v\\
w&z
\end{array}\right) = \;A_k^t \left(\begin{array}{cc}d_1&n\\
n&d_2
\end{array}\right) A_k\;:\; u,v,w,z,k\in\NN\right\}.
$$
It is not difficult to verify that for each $\vn\in N_c$ there
exists $k_\vn\in\NN$ such that
$$
\MMM_\vn = \left\{A_k^t \left(\begin{array}{cc}d_1&n\\
n&d_2
\end{array}\right) A_k\;:\; k\ge k_\vn\right\}.
$$
Finally denote by $\HHH_c$ the set of matrices
$$
\HHH_c:= \{M\in \MMM_\vn\;:\; \vn\in N_c\}.
$$
\begin{proposition}\label{th3}
Let $c<0$. There exists a bijection between
$(\SSS\setminus\{\mathrm{Id}\})\times \HHH_a$ and the set $\FFF_c$
of factorisations of $n^2+c$:
$$
\FFF_a:=\{(n,d_1,d_2)\in \ZZ_{\ge 0}\times \NN^2\;:\; n^2+c =
d_1d_2,\; d_1\le d_2\}.
$$
It is given by:
$$
\phi(M,H):= (n,d_1,d_2),\quad \text{where}\quad \left(\begin{array}{cc}d_1&n\\
n&d_2
\end{array}\right) = M^tHM.
$$
\end{proposition}

The proof of this proposition is analogous to that of
Proposition~\ref{th2} and is left to the interested reader.

Propositions~\ref{th2} and~\ref{th3} provide the way to generate
factorisations of numbers of the form $n^2+c$ for a fixed $c\in
\ZZ$. We provide the schematic description of the algorithm and
leave the details to the reader.

\begin{enumerate}
\item Construct the set $\HHH_c$ for $c>0$ or $N_c$ for $c<0$. This
can be done with help of Lemmata~\ref{lem1} or~\ref{lem2}
respectively.

\item Randomly choose the matrix $H\in \HHH_c$. For $c>0$ this step
is straightforward. For $c\le 0$, one generates the random triple
$\vn\in N_c$ and then chooses $k\ge k_\vn$ to generate $H =
A_k^t\begin{bmatrix}d_1&n\\n&d_2\end{bmatrix}A_k$.

\item Randomly generate $M\in \SSS$ and compute $M^tHM$. The ways to
generate $M$ are discussed in Section~\ref{sec4}.
\end{enumerate}

\bigskip
\noindent Dzmitry Badziahin\\ \noindent The University of Sydney\\
\noindent Camperdown 2006, NSW (Australia)\\
\noindent {\tt dzmitry.badziahin@sydney.edu.au}

\end{document}